\documentclass[letterpaper, 10 pt, conference]{ieeeconf}  % Comment this line out
\IEEEoverridecommandlockouts                              % This command is only
\usepackage{graphics} % for pdf, bitmapped graphics files
\usepackage{epsfig} % for postscript graphics files
\usepackage{mathptmx} % assumes new font selection scheme installed
\usepackage{times} % assumes new font selection scheme installed
\usepackage{amsmath} % assumes amsmath package installed
\usepackage{amssymb}  % assumes amsmath package installed
\DeclareMathOperator*{\argmin}{argmin} 
\usepackage{listings}
\usepackage{xcolor}

\definecolor{codegreen}{rgb}{0,0.6,0}
\definecolor{codegray}{rgb}{0.5,0.5,0.5}
\definecolor{codepurple}{rgb}{0.58,0,0.82}
\definecolor{backcolour}{rgb}{0.95,0.95,0.92}

\lstdefinestyle{mystyle}{
    backgroundcolor=\color{backcolour},   
    commentstyle=\color{codegreen},
    keywordstyle=\color{magenta},
    numberstyle=\tiny\color{codegray},
    stringstyle=\color{codepurple},
    basicstyle=\ttfamily\footnotesize,
    breakatwhitespace=false,         
    breaklines=true,                 
    captionpos=t,                    
    keepspaces=true,                 
    numbers=none,
    frame=none
    numbersep=5pt,                  
    showspaces=false,                
    showstringspaces=false,
    showtabs=false,                  
    tabsize=3
}

\makeatletter
\let\NAT@parse\undefined
\makeatother
\usepackage{hyperref}
\hypersetup{
    colorlinks=true,
    linkcolor=blue,
    filecolor=magenta,      
    urlcolor=blue,
    citecolor=black,
    pdftitle={Overleaf Example},
    }

\title{\LARGE \bf SINDy with Control: A Tutorial}

\author{Urban Fasel, Eurika Kaiser, J. Nathan Kutz, Bingni W. Brunton, and Steven L. Brunton\thanks{This work was supported by the Air Force Office of Scientific Research (AFOSR FA9550-19-1-0386) and the Army Research Office (ARO W911NF-19-1-0045).}%
\thanks{U. Fasel, E. Kaiser, and S. L. Brunton are with the Department of Mechanical Engineering, University of Washington, USA
        {\tt\small ufasel@uw.edu}}%
\thanks{J. N. Kutz is with the Department of Applied Mathematics, University of Washington, USA}%
\thanks{B. W. Brunton is with the Department of Biology, University of Washington, USA}%
\thanks{Code: \url{https://github.com/urban-fasel/SEIR_SINDY_MPC}}%
}

\begin{document}

\maketitle
\thispagestyle{empty}
\pagestyle{empty}

%%%%%%%%%%%%%%%%%%%%%%%%%%%%%%%%%%%%%%%%%%%%%%%%%%%%%%%%%%%%%%%%%%%%%%%%%%%%%%%%
\begin{abstract}
Many dynamical systems of interest are nonlinear, with examples in turbulence, epidemiology, neuroscience, and finance, making them difficult to control using linear approaches. 
Model predictive control (MPC) is a powerful model-based optimization technique that enables the control of such nonlinear systems with constraints. 
However, modern systems often lack computationally tractable models, motivating the use of system identification techniques to learn accurate and efficient models for real-time control. 
In this tutorial article, we review emerging data-driven methods for model discovery and how they are used for nonlinear MPC. 
In particular, we focus on the sparse identification of nonlinear dynamics (SINDy) algorithm and show how it may be used with MPC on an infectious disease control example. 
We compare the performance against MPC based on a linear dynamic mode decomposition (DMD) model. 
Code is provided to run the tutorial examples and may be modified to extend this data-driven control framework to arbitrary nonlinear systems.\\

\end{abstract}

%%%%%%%%%%%%%%%%%%%%%%%%%%%%%%%%%%%%%%%%%%%%%%%%%%%%%%%%%%%%%%%%%%%%%%%%%%%%%%%%
\section{INTRODUCTION}

%\paragraph{Summary/Significance SINDYC}

Modern systems of interest in turbulence, epidemiology, neuroscience, and finance are high-dimensional and nonlinear, and exhibit multiscale phenomena in both space and time~\cite{Kutz2016book}.
Controlling these nonlinear systems remains an important challenge, as traditional linear control approaches are often insufficient.  
There are several control approaches for nonlinear systems, including model predictive control~\cite{mayne2014automatica} and reinforcement learning~\cite{sutton2018reinforcement}. 
These approaches have different performance tradeoffs and requirements, and careful consideration must be given to factors such as the dimensionality of the system, access to an accurate model of the behavior, and the availability of high-quality and abundant data. 
Both approaches are developing rapidly, driven by advanced algorithms in machine learning and optimization, the rise of big data, and improved computational hardware~\cite{brunton2019}. 

Model predictive control (MPC)~\cite{garcia1989model,morari1999model,mayne2014automatica} is a particularly compelling approach that enables control of strongly nonlinear systems with constraints. 
However, MPC relies on efficient models that accurately represent the dynamics of the system, and these models have remained elusive for many disciplines that lack known governing equations. 
Generally, MPC also suffers from the curse of dimensionality, requiring large computational effort and limiting the applicability to low-dimensional problems, often based on locally linear models.
Increasingly, data-driven approaches are providing a hierarchy of models of various fidelity and complexity, which may be used for MPC.  

System identification has a long and rich history in control theory, and it has a close connection with machine learning, as it builds models from data via regression and optimization~\cite{brunton2019}. 
A wide range of data-driven system identification techniques exist in literature, including state-space models from the eigensystem realization algorithm (ERA)~\cite{juang1985} and other subspace identification methods, Volterra
series~\cite{Brockett1976automatica,Boyd1984imajmci}, autoregressive models~\cite{Akaike1969annals} (e.g. ARX, ARMA, NARX and NARMAX~\cite{Billings2013book} models), and neural network (NN) models~\cite{lippmann1987introduction,draeger1995model,aggelogiannaki2008nonlinear,wang2016combined}. 
Nonlinear models based on machine learning, such as NNs, are increasingly used due to advances in computing power, and recently deep learning has been combined with model predictive control~\cite{lenz2015deepmpc} and applied to several complex systems~\cite{baumeister2018deep,morton2018deep,peitz2019koopman,bieker2020deep}. 
Deep reinforcement learning has also been combined with MPC~\cite{Peng.2009,Zhang2016icra}, yielding impressive results in the large-data limit. 
However, machine learning algorithms often lack interpretability and may suffer from overfitting. 

The recent sparse identification of nonlinear dynamics (SINDy) algorithm~\cite{Brunton2016pnas} provides a data-driven model discovery framework, resulting in interpretable models that avoid overfitting, relying on sparsity-promoting optimization to identify parsimonious models from limited data. 
SINDy is closely related to the dynamic mode decomposition (DMD)~\cite{Schmid2010jfm,Rowley2009jfm,Kutz2016book}, which results in a linear model for the evolution of a few key spatiotemporal coherent structures from high-dimensional time-series data. 
DMD was generalized to include external inputs and control~\cite{Proctor2016siads}, and this approach was later applied to develop SINDy with control~\cite{Kaiser2018prsa}.  
Both approaches have been combined with MPC, enabling interpretable models in the low-data limit for real-time control of nonlinear systems. 

\begin{figure*}[ht]
    \begin{center}
    \includegraphics[width=.945\textwidth]{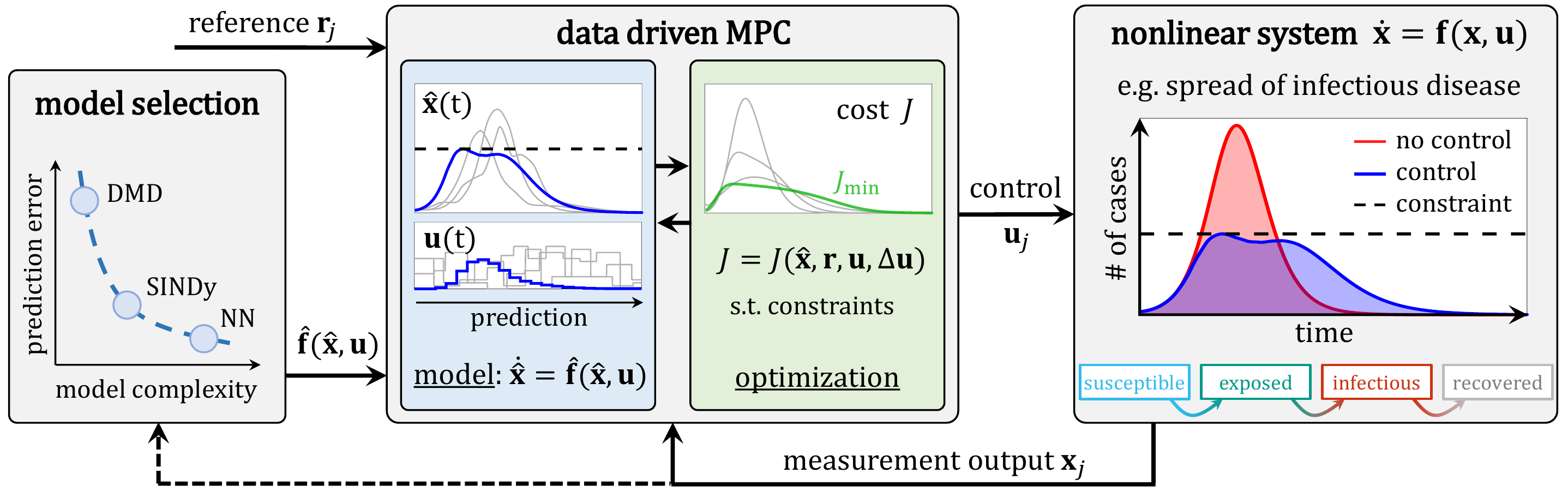}
    \end{center}
    \vspace{-.15in}
    \caption{Schematic of the data-driven MPC framework. Using either linear DMD, sparse nonlinear SINDy, or neural network models for predictive control.}
    \label{fig}
\end{figure*}

In this tutorial article, we review SINDy with control and demonstrate its effectiveness with MPC on an illustrative infectious disease control problem. 
A primary goal of this tutorial is to provide the tools to apply data-driven system identification methods to a model predictive control problem. 
Although we demonstrate the approach with SINDy, many other data-driven modeling techniques may be used, as shown schematically in figure~\ref{fig}. 
Example codes are provided for the tutorial example to promote reproducible research.

%\newpage

%%%%%%%%%%%%%%%%%%%%%%%%%%%%%%%%%%%%%%%%%%%%%%%%%%%%%%%%%%%%%%%%%%%%%%%%%%%%%%%%

\section{DATA-DRIVEN MODELS FOR MPC}

Here we describe data-driven control architectures that combine data-driven model discovery with advanced model-based control strategies. In this tutorial, we use the recent SINDy-MPC architecture~\cite{Kaiser2018prsa} to rapidly identify a low-order model that is used with model predictive control. We consider a nonlinear dynamical system of the form:
\begin{equation}\label{eq1}
    \frac{\mathrm{d}}{\mathrm{dt}}\mathbf{x} = \mathbf{f}(\mathbf{x}, \mathbf{u}), \hspace{10pt} \mathbf{x}(0) = \mathbf{x}_0,
\end{equation}
with state $\mathbf{x}\in\mathbb{R}^n$, control input $\mathbf{u}\in\mathbb{R}^q$ and dynamics $\mathbf{f}(\mathbf{x}, \mathbf{u}): \mathbb{R}^n \times \mathbb{R}^q \rightarrow \mathbb{R}^n$. In this section, we describe SINDy with control and  MPC.  This approach will be used in the next section to provide a step-by-step tutorial and code to compute SINDy-MPC for an infectious disease control problem.

\subsection{Sparse identification of nonlinear dynamics with control}

\begin{figure*}[htbp]
    \includegraphics[width=\textwidth]{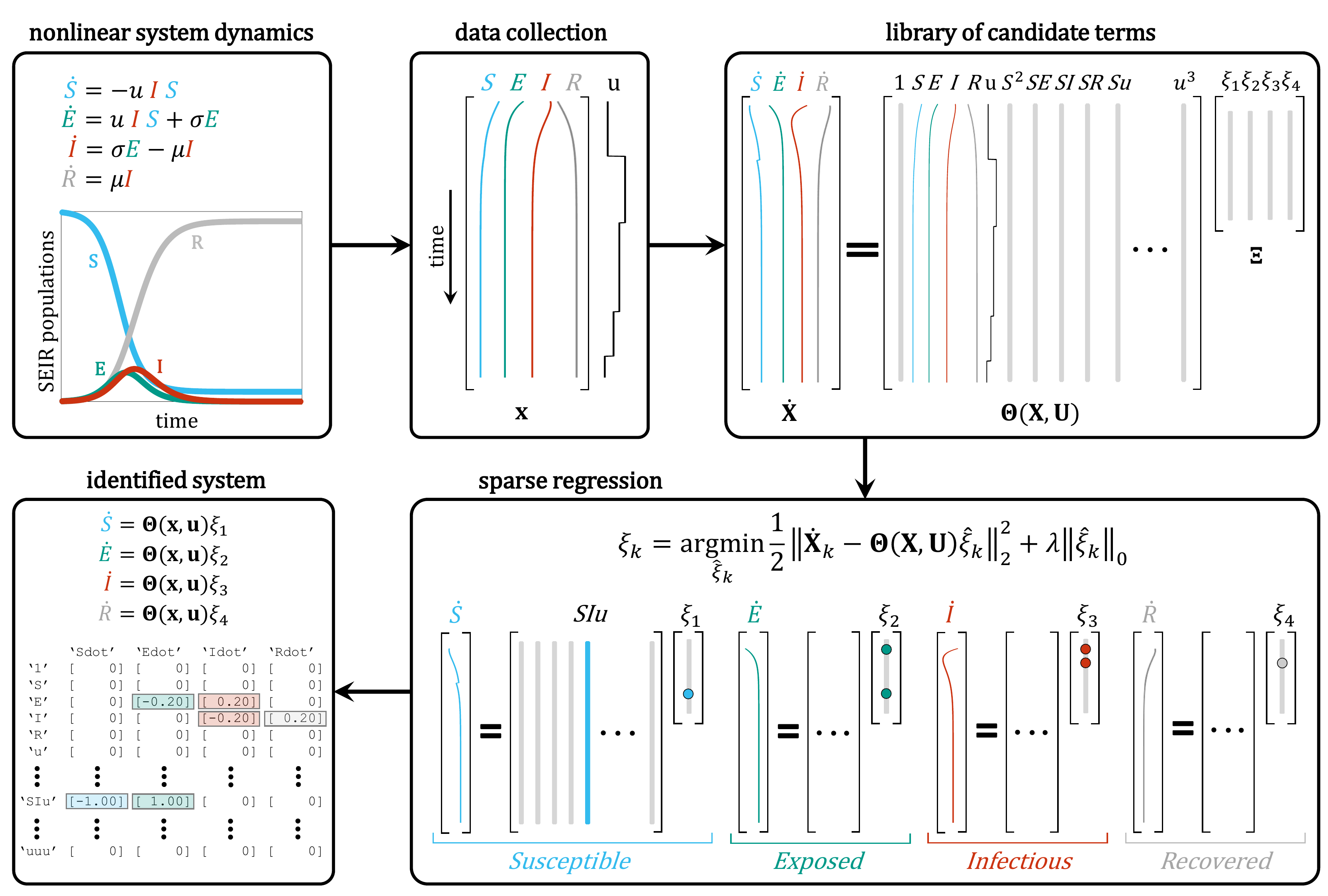}
    \caption{Schematic of the SINDy with control algorithm. Active terms in a library of candidate nonlinearities are selected via sparse regression.}
    \label{fig2}
\end{figure*}

The SINDy~\cite{Brunton2016pnas} and SINDy with control~\cite{Kaiser2018prsa} algorithms identify a sparse nonlinear dynamical system from measurement data, based on the assumption that many systems have relatively few active terms in the dynamics. 
SINDy with control uses sparse regression to identify these few active terms, out of a library $\boldsymbol{\Theta}(\mathbf{x},\mathbf{u})$ of candidate linear and nonlinear model terms in the state $\mathbf{x}$ and actuation $\mathbf{u}$, that are required to approximate the function $\mathbf{f}$ in Eq.~\eqref{eq1}.  
Therefore, sparsity-promoting techniques may be used to find models that automatically balance model complexity with accuracy, resulting in parsimonious models~\cite{Brunton2016pnas}. 
The SINDy with control algorithm is illustrated in figure~\ref{fig2} on a disease model used in the next section. 
To evaluate $\boldsymbol{\Theta}$, we first measure $m$ snapshots of the state $\mathbf{x}$ and the input signal $\mathbf{u}$ in time and arrange these into two matrices\footnote{These matrices are the transposes of the DMD with control matrices~\cite{Proctor2016siads}.}:
\begin{equation}
    \mathbf{X} = \left[\mathbf{x}_1 \hspace{3pt} \mathbf{x}_2 \hspace{3pt} \cdots \hspace{3pt} \mathbf{x}_m\right]^T \hspace{10pt} \mathrm{and} \hspace{10pt} \mathbf{U} = \left[\mathbf{u}_1 \hspace{3pt} \mathbf{u}_2 \hspace{3pt} \cdots \hspace{3pt} \mathbf{u}_m\right]^T.
\end{equation}
After collecting the snapshots, we can evaluate the library of candidate nonlinear functions $\boldsymbol{\Theta}$:
\begin{multline}
    \boldsymbol{\Theta}(\mathbf{X},\mathbf{U}) = [\mathbf{1} \hspace{3pt} \mathbf{X} \hspace{3pt} \mathbf{U} \hspace{3pt} (\mathbf{X} \otimes \mathbf{X}) \hspace{3pt} (\mathbf{X} \otimes \mathbf{U})  \hspace{3pt} (\mathbf{U} \otimes \mathbf{U}) \hspace{3pt} \cdots ],
    %\hspace{3pt} \cdots \\ \hspace{3pt} \sin(\mathbf{X}) \hspace{3pt} \sin(\mathbf{U}) \hspace{3pt} \sin(\mathbf{X} \otimes \mathbf{U}) \hspace{3pt} \cdots ],
\end{multline} 
where $\mathbf{X} \otimes \mathbf{U}$ defines the vector of all product combinations of the components in $\mathbf{x}$ and $\mathbf{u}$. 
This library may include any functions that might describe the data.   
The choice of a suitable library is crucial. The recommended strategy is to start with a basic choice, such as low-order polynomials, and then increase the complexity and order of the library until sparse and accurate models are obtained. 

In addition to evaluating the library, we must compute the time derivatives of the state $\dot{\mathbf{X}} = \left[\dot{\mathbf{x}}_1 \hspace{3pt} \dot{\mathbf{x}}_2 \hspace{3pt} \cdots \hspace{3pt} \dot{\mathbf{x}}_m\right]$, typically by numerical differentiation.  
The system in Eq.~\eqref{eq1} may then be written in terms of these data matrices as:
\begin{equation}
    \dot{\mathbf{X}} = \boldsymbol{\Theta}(\mathbf{X}, \mathbf{U})\boldsymbol{\Xi}.
\end{equation}
Many dynamical systems have relatively few active terms in the governing equations. 
Thus, the coefficients $\boldsymbol{\Xi}$ are mostly sparse and we may employ sparse regression to identify the sparse matrix of coefficients $\boldsymbol{\Xi}$ signifying the fewest nonlinearities in our library that results in a good model fit:
\begin{equation}
    \xi_k = \argmin_{\Hat{\xi}_k} \frac{1}{2} \|\dot{\mathbf{X}}_k-\boldsymbol{\Theta}(\mathbf{X},\mathbf{U})\Hat{\xi}_k\|_2^2 + \lambda\|\Hat{\xi}_k\|_0.
\end{equation}
 $\dot{\mathbf{X}}_k$ is the $k$-th row of $\dot{\mathbf{X}}$, $\xi_k$ is the $k$-th row of $\boldsymbol{\Xi}$, and $\lambda$ is the sparsity-promoting hyperparameter. 
 The term $\| \cdot \|_0$ promotes sparsity in the coefficient vector $\xi_k$, although it is not convex. To approximately solve this optimization, we use sequential thresholded least squares  (STLS)~\cite{Brunton2016pnas}.

 SINDy is closely related to DMD~\cite{Schmid2010jfm,Rowley2009jfm,Kutz2016book}, which extracts spatiotemporal coherent structures from high-dimensional data, along with a linear model for how their amplitudes evolve. 
DMD performs a similar regression to identify a linear discrete-time model $\mathbf{A}$ mapping $\mathbf{X}$ to $\mathbf{X}'$, a matrix with all columns advanced one time step: $\mathbf{X}'=\mathbf{A}\mathbf{X}$. If we formulate SINDy in discrete-time with linear library elements and with $\lambda=0$, SINDy reduces to DMD. It was shown that in many cases, an identified DMD model may be sufficiently accurate to perform control, even for strongly nonlinear dynamical systems~\cite{Kaiser2018prsa}. 
DMD may also provide a useful model until SINDy has collected enough data to accurately characterize the dynamics. 
The \emph{extended DMD} algorithm~\cite{Williams2015jnls}, which seeks linear models in terms of a higher-dimensional state augmented with nonlinear functions of the original variables, has also been used for MPC~\cite{Korda2018automatica}; these DMD-based control approaches are reviewed in~\cite{Brunton2021koopman}.

\subsection{Model predictive control}

This tutorial explores the use of SINDy models for model predictive control. 
MPC is an effective model-based control, which has revolutionized the industrial control landscape~\cite{bryson1975applied,camacho2013model}. MPC enables the control of strongly nonlinear systems with constraints, time delays, non-minimum phase dynamics, and instability. These systems are difficult to control using traditional linear approaches~\cite{garcia1989model,morari1999model,mayne2014automatica}. 

%($\mathbf{y}=\mathbf{x}$ for full state measurements)
In MPC, we compute a control sequence  $\mathbf{u}(\mathbf{x}_j) = \{\mathbf{u}_{j+1},...,\mathbf{u}_{j+m_c}\}$, given the current state estimate or measurement $\mathbf{x}_j$, by a constrained optimization over a receding horizon $T_c = m_c \Delta t$, 
%
% with time step $\Delta t$. % and number of time steps $m_c$. 
with $\Delta t$ the time step of the model and $m_c$ the number of time steps. 
%, which may be different from the sampling time of measurements
%
At each time step, we repeat the optimization, update the control sequence over the control horizon, and apply the first control action $\mathbf{u}_{j+1}$ to the system. The optimal control sequence $\mathbf{u}(\mathbf{x}_j)$ is obtained by minimizing a cost function $J$ over a prediction horizon $T_p = m_p \Delta t \geq T_c$. 
% $m_p$ is the number of time steps of the prediction horizon. The control horizon is generally less than or equal to the prediction horizon, and in case $T_p > T_c$, the input $\mathbf{u}$ is assumed constant thereafter. 
The cost function is: 
%\begin{align*}
%    J = \sum_{k=0}^{m_p-1} \| \hat{\mathbf{x}}_{j+k} - \mathbf{r}_{k} \|{\mathbf{Q}}^2  + \sum_{k=1}^{m_c-1} (\| \mathbf{u}_{j+k} \|{\mathbf{R}_u}^2 + \| \Delta \mathbf{u}_{j+k} \|{\mathbf{R}_{\Delta u}}^2 )\,, 
%\end{align*}
\begin{align*}
    J = \sum_{k=0}^{m_p-1} \| \hat{\mathbf{x}}_{j+k} - \mathbf{r}_{k} \|_{\mathbf{Q}}^2  + \sum_{k=1}^{m_c-1} (\| \mathbf{u}_{j+k} \|_{\mathbf{R}_u}^2 + \| \Delta \mathbf{u}_{j+k} \|_{\mathbf{R}_{\Delta u}}^2 )\,, 
\end{align*}
subject to the discrete-time dynamics and constraints. The cost function $J$ penalizes deviations of the predicted state $\hat{\mathbf{x}}$ from the reference trajectory $\mathbf{r}$, the control expenditure $\mathbf{u}$, and the rate of change of the control signal $\Delta \mathbf{u}$. Each term is weighted by the matrices  $\textbf{Q}$, $\textbf{R}_u$, and $\textbf{R}_{\Delta u}$, respectively. 
To enable this optimization loop to run in real-time, MPC relies on efficient models and high-performance computing. This is particularly challenging for systems where the control must rapidly respond to disturbances on short time-scales, as time delays from sensors, signal transduction, or processing can destroy the robustness of feedback control, putting limitations on the achievable performance~\cite{Doyle2013book}. 

The advantage of MPC lies in simple and intuitive tuning and the ability to control complex phenomena, especially with known constraints and multiple operating conditions.  It also works for systems with time delays and provides the flexibility to formulate and tailor a control objective. 
The major challenge of MPC lies in the development of a suitable model via system identification~\cite{Kaiser2018prsa}. 
%Nonlinear models based on machine learning, such as neural networks, are increasingly used due to advances in computing power, and recently deep reinforcement learning has been combined with MPC~\cite{Peng.2009,Zhang2016icra}, yielding impressive results in the large-data and high-performance computing limit. 
Nonlinear models based on machine learning are increasingly used. However, these techniques often rely on access to massive data sets, have limited ability to generalize, do not readily incorporate known physical constraints, and require expensive and time-consuming computations. 
Instead, Kaiser et al.~\cite{Kaiser2018prsa} showed that simple models obtained via DMD with control~\cite{Proctor2016siads} and SINDy with control~\cite{Kaiser2018prsa} perform nearly as well with MPC on a full nonlinear model, and may be trained in a surprisingly short amount of time.

\section{TUTORIAL ON INFECTIOUS DISEASE}

We now introduce an infectious disease control problem and demonstrate SINDy-MPC to control the spread of the disease. 
We present the main code snippets directly in this tutorial. The complete code to generate the data, identify the SINDy model, and run the MPC can be found on GitHub:
\url{https://github.com/urban-fasel/SEIR_SINDY_MPC}.

\subsection{Intervention strategy for infectious disease}

We first generate the data using an SEIR (susceptible-exposed-infectious-removed) epidemic model~\cite{may1991infectious,castillo1994modeling,chowell2004basic}. 
The SEIR model is a nonlinear system of ordinary differential equations that describes the transition dynamics between four compartments: susceptible $S(t)$ (individuals at risk of
contracting the disease), exposed $E(t)$ (infected individuals but not yet infectious), infectious $I(t)$ (individuals capable of transmitting the disease), and removed $R(t)$ (recovered or dead individuals).
\begin{equation} \label{eqSEIR}
\begin{split}
%\Dot{S}(t) & = -\beta\frac{S(t)I(t)}{N}, \\
%\Dot{E}(t) & = \beta\frac{S(t)I(t)}{N} - k E(t), \\
\Dot{S}(t) & = -\beta S(t)I(t), \\
\Dot{E}(t) & = \beta S(t)I(t) - k E(t), \\
\Dot{I}(t) & = k E(t) -\gamma I(t),  \\
\Dot{R}(t) & = \gamma I(t)
\end{split}
\end{equation}
Susceptible individuals contract the virus at a rate of $\beta I(t)$, with $\beta$ being the transmission rate per day. 
The incubation period $1/k$ determines the time to move from the exposed to the infectious compartment. 
The recovery period $1/\gamma$ determines the time to progress from the infectious to the recovered compartment. 
The basic reproduction number $\mathrm{R_0}=\beta / \gamma$ is an important estimate of the growth of the pandemic. 
$\mathrm{R_0}$ indicates how many infections are generated on average by an infectious individual at the start of a pandemic, when most of the individuals are susceptible. If $\mathrm{R_0} \leq 1$, the infectious cases are declining and the spread of the disease eventually goes to zero. 

To control the spread of an infectious disease, the reproduction number may be reduced through intervention, such as restricted travel, home confinement, or social distancing~\cite{stewart2020control,grimm2021extensions}. 
To mitigate the spread, $\mathrm{R_0}$ needs to be reduced, and to suppress the spread, $\mathrm{R_0}$ needs to be reduced below one. However, controlling $\mathrm{R_0}$ may come at a significant social and economic cost, resulting in a challenging optimization problem to design intervention policies. 
An in-depth discussion of feedback control of infectious disease in the context of COVID-19 was recently reported by Stewart et al.~\cite{stewart2020control}.

\begin{figure*}[t]
    \includegraphics[width=\textwidth]{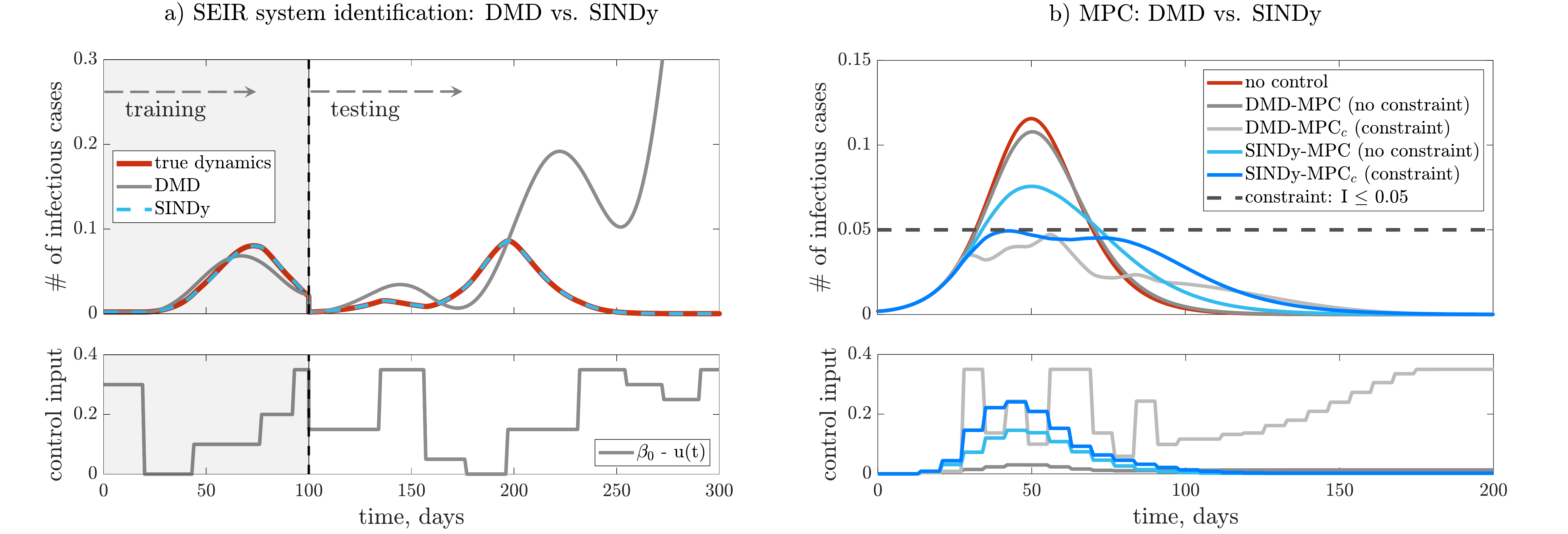}
    \vspace{-.2in}
    \caption{SINDy-MPC vs. DMD-MPC for infectious disease control problem.}
    \label{fig3}
    \vspace{-.2in}
\end{figure*}

\subsection{SINDy-MPC}

In this tutorial, we use SINDy-MPC to control the spread of an arbitrary infectious disease by directly controlling the transmission rate $\beta = \mathrm{u}(t)$. 
Using MPC, we can define constraints on the maximum health care capacity and tailor a control objective considering social cost and economic cost by specifying the cost function weights $\textbf{Q}$, $\textbf{R}_u$, and $\textbf{R}_{\Delta u}$. 
The main challenge of MPC is the development of an accurate and computationally efficient model. 
To identify the nonlinear dynamical system describing the spread of the infectious disease, we use the SINDy with control algorithm, assuming that the dynamical system has relatively few active terms. 
Because these models are sparse by construction, they avoid overfitting, and are more computationally efficient than many other models, so they may be used in real-time and may be identified from relatively small amounts of training data, compared with neural networks. 
The model training and validation, along with the performance of the models for MPC, are shown in figure~\ref{fig3}. 
For comparison, SINDy-based MPC and DMD-based MPC are both presented.  

The state of our system is $\mathbf{x}=[S, E, I, R]$, and we assume that we can control the transmission rate $\beta$ by specific interventions, so the control input is $\mathrm{u}(t) = \beta(t) = \beta_0-\beta_c(t)$. 
The constant parameters of the true SEIR dynamics are set to $\beta_0=0.5$, $\gamma=0.2$, and $k=0.2$. First, we generate the training data with a discrete control input over 100 days. 
The control input is a pseudo-random binary signal (PRBS)~\cite{ljung1999system}, which is a deterministic signal with white-noise-like properties. 
This represents different potential interventions with respective impact on the transmission rate. Code~\ref{code1} generates the data for the forced SEIR dynamics.

\begin{lstlisting}[language=Matlab,caption={Generating SEIR data.},label={code1}]
%% Generate Data
p.b = 0.5; p.g = 0.2; p.k = 0.2; % beta0, gamma, k
n = 4; % Number of states
x0 = [0.996, 0.002, 0.002]; % Initial conditions
tspan = 0:0.1:100;
u = prbsForcing(tspan); % PRBS forcing function
[t,x]=ode45(@(t,x) SEIR(t,x,u,p),tspan,x0);

%% Compute true derivatives
for i=1:length(x)
   dx(i,:) = SEIR(0,x(i,:),u(i),p); 
end

%% SEIR model
function dx = SEIR(t,x,u,p)
   S = x(1); E = x(2); I = x(3); R = x(4);
   dx = [-u*S*I; u*S*I-p.k*E; p.k*E-p.g*I; p.g*I];
end
\end{lstlisting}

After collecting the data, we run the SINDy with control algorithm to identify a sparse nonlinear model. First, we build the library of candidate functions $\boldsymbol{\Theta}$. Here, we use polynomials up to third order. We set the sparsification hyperparameter to $\lambda=0.1$ and run the SINDy algorithm using the sequential threshold least squares approach. Code~\ref{code2} runs SINDy with the functions \texttt{poolData} and \texttt{sparsifyDynamics} that are introduced in~\cite{brunton2019} and can be found on \href{https://github.com/urban-fasel/SEIR_SINDY_MPC}{GitHub}.

\begin{lstlisting}[language=Matlab,caption={Running SINDy with control.},label={code2}]
%% Build library and compute sparse regression
Theta = poolData(x,n,3); % Up to 3rd order polynom
lambda = 0.1; % Sparsification hyperparameter
Xi = sparsifyDynamics(Theta,dx,lambda,n); % STLS
\end{lstlisting}

The control objective is to minimize the number of infected individuals at lowest control expenditure. Additionally, given the maximum health care system capacity (e.g. number of ICU beds), we define a constraint on the number of infected individuals. Here, the reference is $\mathbf{r}=(0,0,0,0)$ and the weight matrices are $\textbf{Q} = \mathrm{diag}(0,0,1,0)$ and $\textbf{R}_u=\textbf{R}_{\Delta u}=0.1$. The control input is limited to $\mathrm{u}=[0.15, 0.5]$, the control is updated once a week and the control and prediction horizon are $m_p=m_c=14\hspace{2pt}\mathrm{days}$. The constraint on the maximum number of infected individuals is set to $5\%$ of the total population. Code~\ref{code3} initializes the MPC problem and runs the SINDy-MPC loop.

\begin{lstlisting}[language=Matlab,caption={Initializing and running SINDy-MPC.},label={code3}]
%% Initialize MPC
pMPC = MPCparams(); % define control parameters
x = x0; uopt = pMPC.uopt0;

%% Run nonlinear MPC with full-state feedback
for i = 1:(pMPC.Duration/pMPC.Ts)   
   % Cost and constraint function
   COST = @(u) CostFCN(u,x,pMPC,uopt(1)); 
   CONS = @(u) ConstraintFCN(u,x,pMPC); 
   % Optimization
   uopt = fmincon(COST,uopt,pMPC,CONS);
   % Apply control and step one timestep forward
   x = rk4u(@SEIR,x,uopt(1),pMPC.Ts,1,[],params);
   xHistory(i,:) = x; 
   uHistory(i) = uopt(1);
end
\end{lstlisting}

In figure~\ref{fig3}, the results of the MPC infectious disease intervention are shown. We use SINDy to identify nonlinear models and compare SINDy-MPC with linear DMD-MPC. On the left of figure~\ref{fig3} in panel a), we show the identification (training, grey background) and prediction (testing, white background) for DMD and SINDy. The control inputs are the same PRBS for SINDy and DMD, shown on the bottom. We see that the SINDy model perfectly predicts the true SEIR dynamics (we only show the infectious population for clarity). The linear DMD model is unstable, diverging after 150 days, and has low predictive performance compared to the SINDy model. We may therefore presume that the DMD model will perform poorly in MPC. 

In panel b) of figure~\ref{fig3}, we show the performance of the SINDy-MPC and compare it with the DMD-MPC. We run both methods with and without constraining the maximum number of infectious individuals (constraint: $\mathrm{I} \leq 0.05$, or $5\%$ of the total population). First, we observe that the SINDy-MPC without the constraint can reduce the number of infectious cases from 11.5\% to 7.5\%. After adding the constraint, the SINDy-$\mathrm{MPC}_c$ can successfully reduce the number of infectious cases below 5\%. We observe that the DMD model under predicts the effect of actuation: the model suggests that further interventions would not reduce the number of cases significantly, and therefore keeps the control and interventions at a lower level. For the constrained case, DMD-$\mathrm{MPC}_c$ is able to reduce the cases below 5\%. The strategy is more uncoordinated and at higher cost compared to the SINDy-$\mathrm{MPC}_c$. However, given the poor predictive accuracy of the DMD model, the DMD-MPC performs surprisingly well, at least over the first third of the spread of the disease. 
We can conclude that DMD models can be of great use, even if their predictive performance is poor. Additionally, these models may be trained with very little data. Therefore, they can be of use in the low-data limit until enough data is collected to identify an accurate SINDy model for control. As concluded in~\cite{Kaiser2018prsa}, SINDy-MPC provides effective and efficient nonlinear control, with DMD as a stopgap until a SINDy model can be identified.

%%%%%%%%%%%%%%%%%%%%%%%%%%%%%%%%%%%%%%%%%%%%%%%%%%%%%%%%%%%%%%%%%%%%%%%%%%%%%%%%

\section{Discussion}

In this tutorial, we explored the use of data-driven model discovery techniques to identify computationally tractable and accurate models of nonlinear systems for model-based control.  
In particular, we demonstrated how SINDy with control can be combined with MPC for infectious disease control. 
We have included example codes throughout to help clarify these concepts.  
Our goal in providing open-source code for this tutorial is to encourage the reader to test the assumptions, explore modifications, and adapt these algorithms to their own nonlinear control problems.

We have made certain assumptions to simplify the SINDy modeling and control procedure. 
We encourage users to implement their own modeling assumptions: changing the numerical differentiation method (assuming we can not measure the derivatives), changing the library functions, the sparse regression algorithm, or investigating different values for the sparsity-promoting hyperparameter $\lambda$ (e.g. using information criteria such as AIC or BIC~\cite{mangan2017model}).
We also encourage the user to investigate different forcing functions and the amount of training data needed to identify models. Forcing the dynamics with single impulses and steps, phase-shifted sum of sinusoids, or other PRBS forcing, will change the condition number of the library $\boldsymbol{\Theta}$, and hence how accurately the model may be identified from limited or noisy data.  
Real-world systems will inevitably have measurement noise and disturbances, which is also important to explore. 
It is also interesting to compare DMD-MPC and SINDy-MPC with MPC based on a neural network. 
Neural networks require significantly more training data and have higher execution time compared to DMD and SINDy. However, they provide a more flexible representation of the dynamics when the model structure varies in state-space~\cite{Kaiser2018prsa}. 
To test this, the code may be modified so the SEIR parameters vary over the course of the spread of the disease. 
We also encourage the user to implement different control strategies (e.g. vaccination and quarantine control~\cite{neilan2010modeling,boujakjian2016modeling}), test different initial conditions, extend the SEIR model to consider other compartments, or completely replace the SEIR dynamics with other nonlinear dynamical systems.

\bibliographystyle{IEEEtran}
\bibliography{references.bib}

\end{document}